\documentclass[11pt]{amsart}

\usepackage[latin2]{inputenc}
\usepackage{t1enc}
\usepackage{amsmath, amsthm, amssymb}
\usepackage{graphicx}
\usepackage{epstopdf}

\newtheorem{Thm}{Theorem}[section]
\newtheorem{Prop}[Thm]{Proposition}
\newtheorem{Lem}[Thm]{Lemma}
\newtheorem{Cor}[Thm]{Corollary}

\theoremstyle{definition}

\newtheorem{Rem}[Thm]{Remark}

\numberwithin{equation}{section}

\newenvironment{Proof}{\rm \trivlist\item[\hskip \labelsep{\bf
Proof.\quad}]}{\hfill\qed\par\medskip\endtrivlist}

\newcommand{\slU}{\mathbf{U}}
\newcommand{\slV}{\mathbf{V}}

\newcommand{\slAb}{\mathbf{Ab}}

\newcommand{\ol}{\overline}

\def\id{\operatorname{id}}

\def\act#1#2{{^{#1}\kern -2pt {#2}}}
\def\acth#1#2{{^{#1}\kern -1pt {#2}}}

\def\graph{\Gamma} %a graph
\def\graphb{\Delta} %another one
\def\connp#1{(S_#1)} %contains no breaking path
\def\freecat#1#2{F_{g#1}(#2)} %free category of paths

\title{On the graph condition regarding the $F$-inverse cover problem}

\author{N\'ora Szak\'acs }
\email{szakacsn@math.u-szeged.hu}
\address{Bolyai Institute\\
University of Szeged\\
Aradi v\'ertan\'uk tere 1\\
H-6720 Szeged, Hungary}

\thanks{This research was supported by the Hungarian National Foundation
for Scientific Research grant no. K104251.}

\begin{document}

\begin{abstract}
In \cite{ASzM}, Auinger and Szendrei have shown that every finite inverse monoid
has an $F$-inverse cover if and only if each finite graph admits a locally finite 
group variety with a certain property. We study this property and prove that
the class of graphs
for which a given group variety has the required property is closed downwards in the minor ordering,
and can therefore be described by forbidden minors. We find these forbidden minors for all
varieties of Abelian groups, thus describing the graphs for which such a group variety satisfies
the above mentioned condition.
\end{abstract}

\maketitle

\section{Introduction}

An {\em inverse monoid} is a monoid $M$ with the
property that for each $a \in M$ there exists a unique element
$a^{-1} \in M$ (the inverse of $a$) such that $a =aa^{-1}a$ and
$a^{-1} = a^{-1}aa^{-1}$. Every inverse monoid may be embedded in a
suitable symmetric inverse monoid $I_V$. Here $I_V$ is the
monoid of all partial injective maps from  $V$ to $V$ (i.e.
bijections between subsets of the set $V$) with respect to the usual
composition of partial maps.

We refer the reader to the books by Lawson \cite{law} or Petrich
\cite{Pet} for the basics of the theory of inverse monoids. In particular,
the {\em natural partial order} on an inverse monoid $M$ is defined
as follows: $a \leq b$ if $a = eb$ for some idempotent $e \in M$.
In $I_V$, this partial order is the one defined by the restriction of partial maps. 
We also recall that inverse monoids, like groups, form a variety of algebras of type $(2,1,0)$, 
and free inverse monoids exist on any set. The free inverse monoid on
the set $X$ is, like the free group, obtained as a factor of the free monoid 
with the involution $^{-1}$, denoted by $FMI(X)$ (see \cite{law} for details).
%, and is denoted by $\FIM (X)$ .

It is well known that each inverse monoid admits a smallest group congruence which is usually denoted by $\sigma$.
An inverse monoid is \emph{$F$-inverse} if each $\sigma$-class has a greatest element with respect to the natural
partial order.

The notion of an $F$-inverse monoid is among the most important ones in the theory of inverse semigroups, for example, free inverse monoids are 
$F$-inverse \cite{law, Pet}. Moreover, they play an important role in the theory of partial actions of groups, see Kellendonk and Lawson
\cite{KelLaw}, and in this context they implicitly occur in Dehornoy \cite{Deh1, Deh2}. In Kaarli and M\'{a}rki \cite{KM}, they occur in the context of universal algebra.
Even in analysis they are useful: see Nica \cite{Nica},
 Khoshkam and Skandalis \cite{KS} and Steinberg \cite{Stein} for their role in the context of $C^\ast$-algebras.

An $F$-inverse monoid $F$ is an \emph{$F$-inverse cover} of the inverse monoid $M$ if there exists an idempotent separating
surjective homomorphism from $F$ to $M$. It is well known that every inverse monoid has an $F$-inverse cover. 
The proof is quite simple and constructive, the inverse cover it yields is an idempotent pure factor of a free inverse monoid, 
and therefore is always infinite.
%A sketch of the proof is to start with
%a surjective homomorphism $\varphi \colon \FIM(X) \to M$, and then factor
%factoring $\FIM(X)$ by the congruence $\theta$ generated by pairs of idempotens with the same image. It can be checked that the factor is $F$-inverse, and
%therefore it is an $F$-inverse cover of $M$ with the mapping $\hat\varphi \colon \FIM(X)/\theta \to M$ induced by $\varphi$.
The question of whether finite inverse monoids admit a finite $F$-inverse cover
was first proposed by Henckell and Rhodes \cite{HR}, and has become
one of the biggest open problems regarding finite inverse semigroups since. 

In \cite{ASzM}, Auinger and Szendrei have translated the $F$-inverse cover problem to the language of group varieties and graphs. They have
proven that the question is equivalent to whether there is, for every finite graph, a locally finite group variety for which a certain condition is satisfied.
These results are summarized in Section 2, and Section 3 contains general observations regarding this condition, including the fact that the class of graphs
for which a given group variety has the required property
is closed downwards in the minor ordering, and can therefore be described by forbidden minors. Section 4 contains our main result, 
which, using forbidden minors, describes the graphs for which there is a variety of Abelian groups satisfying the required condition. 
Not surprisingly, it turns out that these graphs consist of a quite
narrow segment of all finite graphs.

\section{Preliminaries}

We define graphs in this paper to be finite and directed. The set of vertices of a graph $\graph$ is denoted by $V(\graph)$, the set
of edges by $E(\graph)$, and the initial and terminal vertices of an edge $e$ are denoted by $\iota e$ and $\tau e$
respectively. We also say that $e$ is a $(\iota e, \tau e)$-edge. 
We consider $\graph$ as the union of its vertices and edges, the union and intersection of subgraphs of $\graph$ are meant in this sense.
Connectedness of graphs will, however, be
regarded in an undirected sense throughout the paper, that is, we call a digraph \emph{connected} (\emph{two-edge-connected}) if the underlying
undirected graph is connected (two-edge-connected). Recall that an undirected graph is called \emph{two-edge-connected} if it is connected
and remains connected whenever an edge is removed.

A \emph{path} in $\graph$ is a sequence $e_1 \cdots e_n$ of consecutive edges which, by definition,
means that $\tau e_i =\iota e_{i+1}$, $i=1, \ldots, n-1$, or an empty path around an arbitrary vertex. 
There is an evident notion of inital and terminal vertices of paths, also denoted by $\iota$ and $\tau$ respectively:
if $p=e_1 \cdots e_n$ for some $n \in \mathbb N$, then $\iota p=\iota e_1, \tau p= \tau e_n$, and if $p$ is an empty path around $v$, then $\iota p=\tau p=v$.
We do, however, also need to consider paths in a more general,
undirected sense. We therefore introduce the graph $\ol\graph$, for which $V(\ol\graph)=V(\graph)$ and
$E(\ol\graph)=E(\graph) \cup (E(\graph))^{-1}$, where $E(\graph) \cap (E(\graph))^{-1}=\emptyset$, and for an edge $e \in E(\graph)$, $e^{-1}$ is a $(\tau e, \iota e)$-edge.
We will often consider paths in $\ol\graph$. 
In this context, we can also define the inverse of a path $p$ in $\ol\graph$, denoted by $p^{-1}$, to be the $(\tau p, \iota p)$-path traversing the edges of $p$ in the opposite direction.
Note that a path in $\ol\graph$ can be regarded as a word in the free monoid
$FMI(E(\graph))$ with involution, where, of course, the empty paths correspond to the empty word,
and the word corresponding to the path $p$ is of course the inverse of the words corresponding to $p^{-1}.$
For a path $p$ in $\ol\graph$, the \emph{graph spanned by $p$}, denoted by $\langle p \rangle$,
is the subgraph of $\graph$ consisting of the vertices $p$ traverses and edges $e \in E(\graph)$ for which $e$ or $e^{-1}$ occurs in $p$.

An \emph{inverse category} is a category in which each arrow $x$ admits a unique arrow $x^{-1}$ satisfying $x=xx^{-1}x$ and $x^{-1}=x^{-1}xx^{-1}$.
The natural partial order on inverse categories is defined the same way as for inverse monoids, that is, $x \leq y$ if $x=ey$ for some
idempotent arrow $e$.

We summarize some definitions and results of \cite{ASzM} necessary to formulate our results.
Let $\slV$ be a variety of inverse monoids and let $X$ be an alphabet. 
%We denote by $F_\slV(X)$ the relatively free inverse monoid in $\slV$ on $X$.
%For words $u,v \in (X\cup X^{-1})^\ast$, we denote by $[u]_\slV$ the value of $u$ in $F_\slV(X)$, and put $u \equiv_\slV v$ if $[u]_\slV=[v]_\slV$.
For words $u, v \in FMI(X)$, we put $u \equiv_\slV v$ if the identity $u=v$ holds in $\slV$. It is well known that $\equiv_\slV$ is a fully invariant
congruence on $FMI(X)$, and $F_\slV(X)=FMI(X)/\equiv_\slV$ is the relatively free inverse monoid in $\slV$ on $X$. We denote by $[u]_\slV$ the $\equiv_\slV$-class
of $u$, that is, the value of $u$ in $F_\slV(X)$.

Let $\graph$ be a graph and $\slV$ be a variety of inverse monoids. By $\freecat \slV \graph$, we denote the \emph{free $g\slV$-category} on $\graph$:
its set of vertices is $V(\graph)$, its set of $(i,j)$-arrows is
$$\freecat\slV\graph(i,j)=\{(i, [p]_\slV, j) : p \hbox{ is an } (i,j)\hbox{-path in } \ol\graph \},$$
and the product of two consecutive arrows is defined by
$$(i, [p]_\slV, j) (j, [p]_\slV, k)=(i, [p]_\slV[q]_\slV, k)=(i, [pq]_\slV, k).$$
The inverse of an arrow is given by
$$(i, [p]_\slV, j)^{-1}=(j, [p]^{-1}_\slV, i)=(j, [p^{-1}]_\slV, i).$$

An important case is $\slV=\mathbf{Sl}$, the variety of semilattices, in which case $[p]_{\mathbf{Sl}}$ can be identified with the subgraph $\langle p \rangle$ of $\graph$ spanned by $p$.
The $(i,j)$-arrows of $\freecat{\mathbf{Sl}}\graph$ are therefore precisely the triples $(i, \graphb,j)$, where $\graphb$ is a
connected subgraph containing $i$ and $j$.
The natural partial order on $\freecat{\mathbf{Sl}}\graph$ is conveniently described as
$$(i, \graphb_1, j) \leq (k, \graphb_2, l) \hbox{ if and only if } i=k, j=l \hbox{ and } \graphb_1 \supseteq \graphb_2.$$

A \emph{dual premorphism} $\psi \colon C \to D$ between inverse categories is a graph homomorphism satisfying $(x\psi)^{-1}=x^{-1}\psi$ and $(xy)\psi \geq x\psi \cdot y\psi$.
According to \cite{ASzM}, every finite inverse monoid admits a finite $F$-inverse cover if and only if, for every (finite) connected graph $\graph$, 
there exist a locally finite group variety $\slU$ and a dual premorphism 
$\psi \colon \freecat\slU\graph \to \freecat{\mathbf{Sl}}\graph$ with $\psi| \graph=\id_\graph$.

Now fix a connected graph $\graph$ and a group variety $\slU$. We assign to each arrow $x$ of $\freecat\slU\graph$ two sequences of finite subgraphs of $\graph$ as follows: let
\begin{equation}
\label{C_0}
C_0(x)=\bigcap\{\langle p \rangle  : (\iota p, [p]_\slU, \tau p)=x\},
\end{equation}
and let $P_0(x)$ be the connected component of $C_0(x)$ containing $\iota x$. If $C_n(x), P_n(x)$ are already defined for all $x$, then put
$$C_{n+1}(x)=\bigcap\{P_n(x_1)\cup \cdots \cup P_n(x_k) : k\in \mathbb N, x_1 \cdots x_k=x\},$$
and again, $P_{n+1}(x)$ is the connected component of $C_{n+1}(x)$ containing $\iota x$. 

It is easy to see that 
$$C_0(x) \supseteq P_0(x) \supseteq \cdots \supseteq C_n(x) \supseteq P_n(x) \supseteq C_{n+1}(x) \supseteq P_{n+1}(x) \supseteq \cdots$$
for all $x$ and $n$. We define $P(x)$ to be $\bigcap_{n=0}^\infty P_n(x)$, which is a connected subgraph of $\graph$ containing $\iota x$.
According to 
\cite[Lemma 3.1]{ASzM}, 
there exists a dual premorphism $\psi \colon \freecat\slU\graph \to \freecat{\mathbf{Sl}}\graph$
with $\psi|\graph=\id_\graph$ if and only if $\tau x \in P(x)$ for all $x$, and in this case, the assignment
$x \mapsto (\iota x, P(x), \tau x)$ gives such a dual premorphism.
If $\tau x \notin P(x)$ for some $x=(\iota p, [p]_\slU, \tau p)$, then we call $p$ a \emph{breaking path} over $\slU$.

In \cite{ASzM}, $C_0(x)$ is incorrectly defined to be the graph spanned by the $\slU$-content of $x$ together with $\iota x$. From
the proof of \cite[Lemma 3.1]{ASzM} (see the inclusion $\mu(x\psi)\subseteq C_0(x)$), it is clear that the definition of $C_0(x)$ needed is the one in (\ref{C_0}).  
The following proposition states that in the cases crucial for the main result 
\cite[Theorem 5.1]{ASzM},
i.e., where $\graph$ is the Cayley graph of a finite group, these two definitions are equivalent in the sense that $P_0(x)$, and so the sequence $P_n(x)$
does not depend on which definition we use.
For our later convenience, let $\hat C_0(x)$ denote the graph which is the union of the $\slU$-content of $x$ and $\iota x$.

\begin{Lem}
\label{biconn}
If $\graph$ is two-edge-connected, then for any arrow $x$ of $\freecat\slU\graph$, the subgraphs $C_0(x)$ and $\hat C_0(x)$ can only differ in isolated vertices
(distinct from $\iota x$ and $\tau x$).
\end{Lem}

\begin{Proof}
Let $x$ be an arrow of $\freecat\slU\graph$.
It is clear that $\hat C_0(x) \subseteq C_0(x)$. For the converse, put $x=(\iota p, [p]_\slU, \tau p)$, and suppose $e$ is an edge of $\langle p \rangle$
such that $e \notin \hat C_0(x)$. Let $s_e$ be a $(\iota e,\tau e)$-path in $\ol\graph$ not containing $e$ --- such a path exists since $\graph$ is two-edge-connected.
Let $p_{e \to s_e}$ be the path obtained from $p$ by replacing all occurences of $e$ by $s_e$. Then $p \equiv_\slU p_{e \to s_e}$, and $e \notin \langle p_{e \to s_e} \rangle$,
hence $e \notin C_0(x)$, which completes the proof.
\end{Proof}

\begin{Rem}
We remark that the condition of $\graph$ being two-edge-connected is necessary in Lemma \ref{biconn}, that is,
when $\graph$ is not two-edge-connected, the subgraphs $C_0(x)$ and $\hat C_0(x)$ can in fact be different. Put, for example,
$\slU=\slAb$, the variety of Abelian groups, and let $e$ be an edge of $\graph$ for which $\graph\backslash \{e\}$ is disconnected. Let $p=ese^{-1}$ be a path in $\ol\graph$,
where $s \not\equiv_\slAb 1$ and $e, e^{-1}$ do not occur in $s$. Then the subgraph spanned by the $\slAb$-content of $p$ does not contain $e$, whereas any path $p'$ which is coterminal with and $\slAb$-equivalent
to $p$ must contain the edge $e$, as there is no other $(\iota e ,\tau e)$-path in $\ol\graph$. 
\end{Rem}

For a group variety $\slU$, we say that a graph $\graph$ satisfies property $\connp\slU$,
or $\graph$ is $\connp\slU$ for short, if $\tau x \in P(x)$ holds for any arrow $x$ of $\freecat\slU\graph$. 
By \cite{ASzM}, each finite inverse monoid has a finite $F$-inverse cover if and only if each finite connected graph is $\connp\slU$ for some locally finite group variety $\slU$.
This property $\connp\slU$ for finite connected graphs is our topic for the remaining part of the paper.

We recall that by \cite[Lemmas 4.1 and 4.2]{ASzM}, the following holds.

\begin{Lem}
\label{subgraph}
If a graph $\graph$ is $\connp\slU$ for some group variety $\slU$, then so is any redirection of $\graph$, and any subgraph of $\graph$.
\end{Lem}

However, we remark that the lemma following these observations in \cite{ASzM}, namely Lemma 4.3 is false. It states that if a simple graph $\graph$ is $\connp\slU$, then so is any graph obtained from
$\graph$ by adding parallel edges (where both ``simple'' and ``parallel'' are meant in the undirected sense). Our main result Theorem \ref{mainabel} yields counterexamples.

\begin{Lem}
\label{su-sv}
If\/ $\slU$ and $\slV$ are group varieties for which $\slU \subseteq \slV$, then $\connp \slU$ implies $\connp \slV$.
\end{Lem}

\begin{Proof}
Suppose $\graph$ is $\connp \slU$, let $p$ be any path in $\ol\graph$. Put $x^\slU=(\iota p, [p]_\slU, \tau p) \in \freecat\slU\graph$, and similarly 
let $x^\slV=(\iota p, [p]_\slV, \tau p) \in \freecat\slV\graph$. Since 
$\slU \subseteq \slV$, we have $C_0(x^\slU) \subseteq C_0(x^\slV)$. Also, since $q \equiv_\slV q_1 \cdots q_n$ implies $q \equiv_\slU q_1 \cdots q_n$, 
we obtain $P_n(x^\slU) \subseteq P_n(x^\slV)$ by induction. Since $\tau p \in P_n(x^\slU)$ by assumption, this yields $\tau p \in P_n(x^\slV)$, that is, $\graph$ is $\connp \slV$.
\end{Proof}

%
%
%
%%
%%%
%%%% NEW SECTION !!!
%%%
%%
%
%
%

\section{Forbidden minors}

In this section, we prove that, given a group variety $\slU$, the class of graphs satisfying $\connp\slU$ can be described by forbidden minors.

Let $\graph$ be a graph and let $e$ be a $(u,v)$-edge of $\graph$ such that $u \neq v$. The operation which removes $e$ and simultaneously merges $u$ and $v$
to one vertex is called \emph{edge-contraction}. We call $\graphb$ a \emph{minor} of $\graph$ if it can be obtained from $\graph$ by edge-contraction, 
omitting vertices and edges, and redirecting edges.

\begin{Prop}
Suppose $\graph$ and $\graphb$ are connected graphs such that $\graphb$ is a minor of $\graph$. Then, if $\graphb$ is non-$\connp \slU$, 
so is $\graph$.
\end{Prop}

\begin{Proof}
By Lemma \ref{subgraph}, adding edges and vertices to, or redirecting some edges of a graph does not change the fact that it is non-$\connp \slU$. Therefore 
let us suppose that $\graphb$ is obtained from $\graph$ by contracting an edge $e$ for which $\iota e \neq \tau e$.
Let $x_1, \ldots, x_n$ be the edges of $\graph$ having $\iota e$
as their terminal vertex. 
For a path $p$ in $\ol \graphb$, let $p_{+e}$ denote the path in $\ol \graph$ obtained by replacing all occurrences of $x_j$  $(j=1, \ldots, n)$ by $x_je$ (and
all occurences of $x_j^{-1}$ by $e^{-1}x_j^{-1}$).
Similarly, for a subgraph $\graphb'$ of $\graphb$, let $\graphb'_{+e}$ denote the subgraph of $\graph$ obtained from $\graphb'$ by taking its preimage
under the edge-contraction containing the edge $e$ if $\graphb'$ contains some $x_j$ $(j=1, \ldots, n)$, and its preimage without $e$ otherwise. Obviously, we have
$\langle p_{+e} \rangle = \langle p \rangle_{+e}$ for any path $p$ in $\ol \graphb$.

%Since for any paths $q, q_1, \ldots, q_k$ in $\ol \graphb$, we have $q \equiv_\slAb q_1 \cdots q_n$ if and only if $q_{+e} \equiv_\slAb (q_1)_{+e} \cdots (q_n)_{+e}$,
%we obtain $(C_0(p))_{+e}=C_0(p_{+e})$ for any path $p$, furthermore $(P_n(p))_{+e}=P_n(p_{+e})$ holds by induction. 
Note that if $p$ is a path in $\ol\graphb$ traversing the edges $f_1, \ldots, f_k$, then $p_{+e}$, considered as a word in $FMI(\{e, f_1, \ldots, f_k\})$,
is obtained from the word $p$ by substituting $(x_je)$ for $x_j$ $(j=1,\ldots, n)$, 
and leaving the other edges unchanged. Putting $x=(\iota p, [p]_\slU, \tau p)$ and $x_{+e}=(\iota p_{+e}, [p_{+e}]_\slU, \tau p_{+e})$, this implies $(C_0(x))_{+e} \supseteq C_0(x_{+e})$
for any path $p$ is $\ol\graphb$. Moreover, we also see that, for any paths $q, q_1, \ldots, q_k$ in $\ol \graphb$, 
we have $q \equiv_\slU q_1 \cdots q_n$ if and only if $q_{+e} \equiv_\slU (q_1)_{+e} \cdots (q_n)_{+e}$. 
Using that for any subgraph $\graphb' \subseteq \graphb$,
the connected components of $\graphb'$ and $\graphb'_{+e}$ are in one-one correspondence, an induction shows that $(P_n(x))_{+e}\supseteq P_n(x_{+e})$ for every $n$.
In particular, $P_n(x)$ contains $\tau p$ if and only if 
$(P_n(x))_{+e}$ contains $\tau p_{+e}$. Therefore if $p$ is a breaking path in $\graphb$ over $\slU$, then $\tau p_{+e} \notin (P_n(x))_{+e}$ 
and hence $\tau p_{+e} \notin P_n(x_{+e})$, that is, $p_{+e}$ is a breaking path in $\graph$ over $\slU$, which proves our 
statement.
\end{Proof}

By the previous proposition, the class of all graphs containing a breaking path over $\mathbb \slU$ (that is, of all non-$\connp \slU$ graphs) is closed upwards in the minor ordering,
hence, it is determined by its minimal elements. According to the theorem of Robertson and Seymour \cite{RobSer}, there is no infinite anti-chain in the minor
ordering, that is, the set of minimal non-$\connp \slU$ graphs must be finite.

These observations are summarized in the following theorem:

\begin{Thm}
For any group variety $\slU$, there exist a finite set of connected graphs $\graph_1, \ldots, \graph_n$ such that 
the graphs containing a breaking path over $\mathbb \slU$ are exactly those having one of $\graph_1, \ldots, \graph_n$ as a minor.
\end{Thm}

By Lemma \ref{su-sv}, if $\slU$ and $\slV$ are group varieties with $\slU \subseteq \slV$, the forbidden minors for $\slU$ are smaller (in the minor ordering) then the ones for $\slV$.

The next statement contains simple observations regarding the nature of forbidden minors.

\begin{Prop}
\label{mingraphs}
For any group variety $\slU$, the set of minimal non-$\connp\slU$ graphs are two-edge-connected graphs without loops.
\end{Prop}

\begin{Proof}
We show that if $\graph$ is a non-$\connp\slU$ graph which has loops or is not two-edge-connected, then there exists a graph below $\graph$ in the minor
ordering which is also non-$\connp\slU$. Indeed, suppose that $\graph$ has a loop $e$, and take $\graph\backslash \{e\}$. For a path $p$ in $\ol\graph$, let 
$p_{-e}$ denote the corresponding path in $\ol\graph$ obtained by omitting all occurences of $e$, and for an arrow $x=(\iota p, [p]_\slU, \tau p)
\in \freecat\slU\graph$, put $x_{-e}=(\iota p, [p_{-e}]_\slU, \tau p) \in \freecat\slU{\graph\backslash \{e\}}$. Then 
it is easy to see by induction
that $C_n(x_{-e})\subseteq C_n(x)\backslash \{e\}$ and $P_n(x_{-e})\subseteq P_n(x)\backslash \{e\}$ for every $x$ and $n$, and hence $\tau p \in P_n(x_{-e})$ 
implies $\tau p \in P_n(x)\backslash \{e\}$.

Now suppose $\graph$ is not two-edge-connected, that is, there is a $(u,v)$-edge $e$ of $\graph$ for which $\graph\backslash \{e\}$ is disconnected. Then let $\graph_{u=v}$ denote the
graph which we obtain from $\graph$ by contracting $e$. For a path $p$ in $\ol\graph$, 
let $p_{u=v}$ denote the path in $\ol{\graph_{u=v}}$ which we obtain by omitting all occurrences of $e$ from $p$, and for an arrow $x=(\iota p, [p]_\slU, \tau p)
\in \freecat\slU\graph$, put $x_{u=v}=(\iota p_{u=v}, [p_{u=v}]_\slU, \tau p_{u=v}) \in \freecat\slU{\graph_{u=v}}$.
Observe that for coterminal paths $s,t$ in $\ol\graph$, $s \equiv_\slU t$ implies $s_{u=v}\equiv_\slU t_{u=v}$. This, by induction yields
$C_n(x_{u=v}) \subseteq C_n(x)_{u=v}$ and $P_n(x_{u=v}) \subseteq P_n(x)_{u=v}$ for all $n$, and hence $\tau p \in P_n(x_{u=v})$ implies $\tau p \in P_n(x)_{u=v}$.
\end{Proof}

%
%
%
%%
%%%
%%%% NEW SECTION !!!
%%%
%%
%
%
%

\section{Main result}

In this section, we describe the forbidden minors (in the sense of the previous section) for all non-trivial varieties of Abelian groups.
Denote by $\slAb$ the variety of all Abelian groups.

\begin{center}
\begin{figure}[h]

\includegraphics[width=0.85\linewidth]{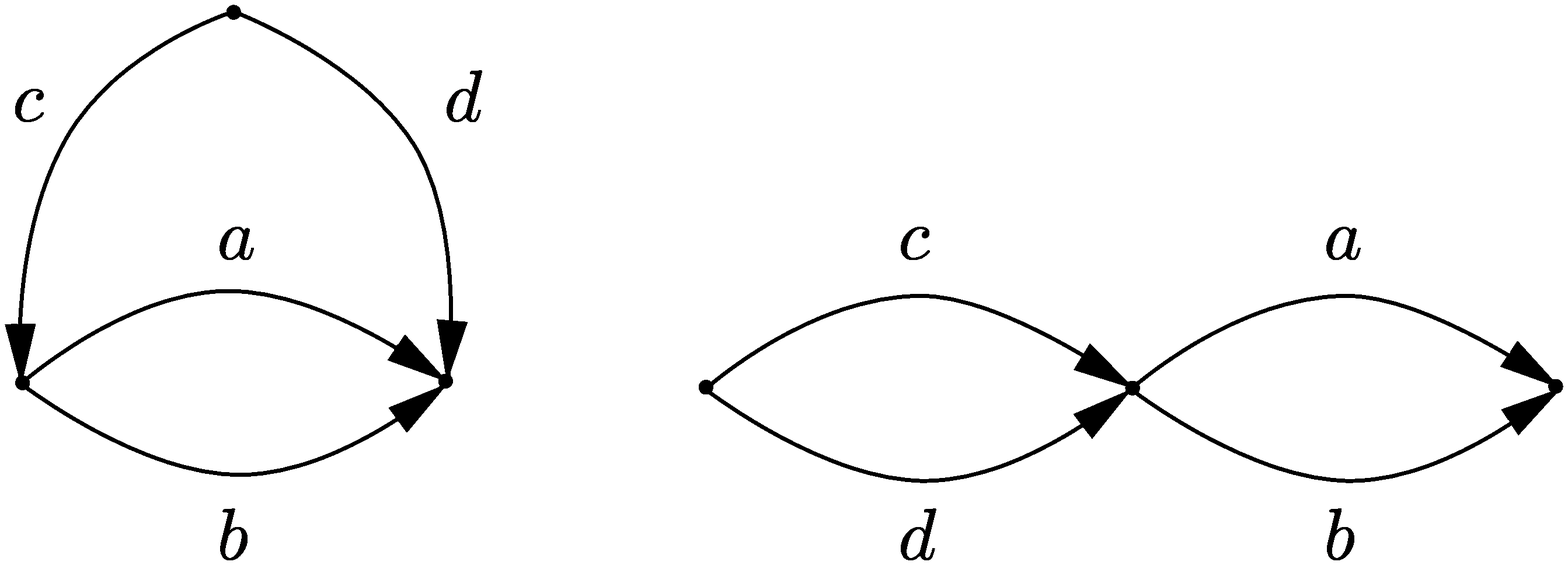}
\caption{}
\label{fig:figure}

\end{figure}
\end{center}

\begin{Thm}
\label{mainabel}
A connected graph contains a breaking path over $\slAb$ if and only if its minors contain at least one of the graphs in Figure \ref{fig:figure}.
\end{Thm}

\begin{Proof}
First, suppose $\graph$ is a finite graph which does not have either graph 
in Figure \ref{fig:figure}
as a minor. Then $\graph$ is either a cycle of length $n$ for some $n \in \mathbb N_0$ with possibly some trees and loops attached, or a graph with at most $2$ vertices. 
According to Proposition \ref{mingraphs}, $\graph$ contains a breaking path if and only if its greatest two-edge-connected minor does, which, in the formes case is the cycle $\graph_n$ of length $n$, and in the latter case is a two-edge-connected graph on at most $2$ vertices. 
It is easy to see that both in cycles $\graph_n$ or graphs on at most $2$ vertices, for any path $p$,
the $\slAb$-content $\hat C_0(x)$ with $x=(\iota p, [p]_\slAb, \tau p)$ is connected, therefore by Lemma \ref{biconn}, these graphs do not contain a breaking path
over $\slAb$.

For the converse part, we prove that both graphs 
in Figure \ref{fig:figure}
contain a breaking path over $\slAb$ --- namely, the path $a$. 
For brevity, denote $\iota a, \tau a$ and $\iota c$ by $u, v$ and $w$ respectively, and put $x=(u, [a]_\slAb, v)$.
Since both graphs are two-edge-connected, 
Lemma \ref{biconn} implies that $C_0(x)$ and the $\slAb$-content $\hat C_0(x)=\langle a \rangle$ are (almost) the same, that is,
$P_0(x)=\langle a \rangle$ in both cases.
Now put $x_1=( u, [c^{-1}]_\slAb, w),\ x_2=(w,[cab^{-1}c^{-1}]_\slAb ,w),\ x_3=(w,[cb]_\slAb,v)$, and note that $x=x_1x_2x_3$, that is,
$C_1(x) \subseteq P_0(x) \cap (P_0(x_1) \cup P_0(x_2) \cup P_0(x_3))$.
Again, using $\hat C_0$ and Lemma \ref{biconn}, we obtain that 
$\hat C_0(x_1)=\langle c \rangle$, 
$\hat C_0(x_2)=\{w\} \cup \langle ab^{-1} \rangle$,
$\hat C_0(x_3)=\langle cb \rangle$, and so
$P_0(x_1) \cup P_0(x_2) \cup P_0(x_3)=\langle c \rangle \cup \{w\} \cup \langle cb \rangle=\langle cb \rangle$ for both graphs in Figure \ref{fig:figure}.
Therefore $C_1(x)\subseteq \langle a \rangle \cap \langle cb \rangle=\{u,v\}$ and so $v \notin P_1(x) \subseteq \{u\}$. Hence
$a$ is, indeed, a breaking path over $\slAb$ in both graphs.
\end{Proof}

\begin{Cor}
For any non-trivial variety $\slU$ of Abelian groups, a connected graph contains a breaking path over $\slU$ if and only if its minors contain at least one of 
the graphs in Figure \ref{fig:figure}.
\end{Cor}

\begin{Proof}
The statement is proven in Theorem \ref{mainabel} if $\slU=\slAb$. Now let $\slU$ be a proper subvariety of $\slAb$. 
Then $\slU$ is the variety of Abelian groups of exponent $n$ for some positive integer $n \geq 2$. By Lemma \ref{su-sv}, the forbidden minors for $\slU$ must be minors of 
one of the forbidden minors of $\slAb$, that is, by Proposition \ref{mingraphs}, they are either the same, or the only forbidden minor is the cycle $\graph_2$ of length two. However, it is clear that $\graph_2$ contains
no breaking path over $\slU$ for the same reason as in the case of $\slAb$, which proves our statement.
\end{Proof}

\begin{Rem}
For the variety $\mathbf 1$ of trivial groups, a finite connected graph is $\connp{\mathbf 1}$ if and only if it is a tree with some loops attached. That is, even the smallest two-edge-connected graph
in the minor ordering, the cycle of length two contains a breaking path over $\mathbf 1$.
\end{Rem}

\section*{Acknowledgement}

I would like to thank my supervisor, M\'aria B. Szendrei for introducing me to the problem and giving me helpful advice during my work.

\end{document}